\documentclass[graybox]{svmult}

\usepackage{type1cm}        
\usepackage{makeidx}         
\usepackage{graphicx}        
\usepackage{multicol}        
\usepackage[bottom]{footmisc}

\usepackage{newtxtext}       %
\usepackage{newtxmath}       

\makeindex             

\usepackage{url}
\usepackage{comment}

\allowdisplaybreaks[1]

\renewcommand*{\bar}{\overline}
\newcommand{\gfp}[1]{\Gamma_p{\left({#1}\right)}}
\newcommand{\biggfp}[1]{\Gamma_p{\bigl({#1}\bigr)}}

\makeatletter
\def\imod#1{\allowbreak\mkern5mu({\operator@font mod}\,\,#1)}
\makeatother

\makeatletter
\@namedef{subjclassname@2020}{%
  \textup{2020} Mathematics Subject Classification}
\makeatother

\begin{document}

\motto{In honor of Lance Littlejohn on his 70th birthday.}

\title*{Hypergeometric Functions over Finite Fields and Modular Forms: A Survey and New Conjectures}
\titlerunning{Hypergeometric Functions over Finite Fields and Modular Forms}
\author{Madeline Locus Dawsey and Dermot M\lowercase{c}Carthy}
\institute{Madeline Locus Dawsey \at Department of Mathematics, The University of Texas at Tyler, Tyler, TX 75799, USA,\\ \email{mdawsey@uttyler.edu}
\and Dermot M\lowercase{c}Carthy \at Department of Mathematics \& Statistics, Texas Tech University, Lubbock, TX 79410-1042, USA,\\ \email{dermot.mccarthy@ttu.edu}\\
\thanks{The second author was supported by a grant from the Simons Foundation (\#353329, Dermot McCarthy).}}
%

\maketitle

\abstract*{Hypergeometric functions over finite fields were introduced by Greene in the 1980s as a finite field analogue of classical hypergeometric series. These functions, and their generalizations, naturally lend themselves to, and have been widely used in, character sum evaluations and counting points on algebraic varieties. More interestingly, perhaps, are their links to Fourier coefficients of modular forms. In this paper, we outline the main results in this area and also conjecture 13 new relations.}

\abstract{Hypergeometric functions over finite fields were introduced by Greene in the 1980s as a finite field analogue of classical hypergeometric series. These functions, and their generalizations, naturally lend themselves to, and have been widely used in, character sum evaluations and counting points on algebraic varieties. More interestingly, perhaps, are their links to Fourier coefficients of modular forms. In this paper, we outline the main results in this area and also conjecture 13 new relations.}


\section{Introduction}\label{sec_Intro}

Hypergeometric functions over finite fields were introduced by Greene \cite{G} as a finite field analogue of classical hypergeometric series. 
Much of Greene's early work on these functions focused on developing transformation and summation formulas which mirror those of the classical series. These functions have a nice character sum representation and so the transformation and summation formulas can be interpreted as relations to simplify and evaluate complex character sums. Using character sums to count points on certain algebraic varieties over finite fields is well-established and, consequently, hypergeometric functions over finite fields naturally lend themselves to this endeavor.

Via modularity results, we then find relations between hypergeometric functions over finite fields and the Fourier coefficients of modular forms. As we will see, these relations are striking in their simplicity. While the modularity theorem and connections between hypergeometric functions over finite fields and elliptic curves yield infinitely many such relations with weight two newforms, relations in higher weights are rare. The main source of relations between hypergeometric functions over finite fields and the Fourier coefficients of modular forms of weight greater than two is the supercongruence conjectures of Rodriguez Villegas which yields 14 such relations. In a recent paper studying generalized Paley graphs, the authors discovered evidence for two new such relations. Since then we have conducted a more extensive search where we have found a further 13 possible relations. The main purpose of this paper is to present the details of these conjectural relations. For context, we also outline the other main results in this area.

This paper is organized as follows. In Section \ref{sec_Prelim}, we define hypergeometric functions over finite fields and also a $p$-adic extension.  We then outline the main results linking hypergeometric functions over finite fields and Fourier coefficients of modular forms, categorized by weight, with weight two in Section \ref{sec_wt2} and higher weights in Section \ref{sec_higherwt}. In Section \ref{sec_Trace}, we discuss the Eichler--Selberg trace formula which is one of the main tools used to prove such results. Finally, in Section \ref{sec_NewConjs} we describe our 13 new conjectural relations.


\section{Preliminaries}\label{sec_Prelim}

While hypergeometric functions over finite fields were originally defined by Greene \cite{G}, in this paper we will use a normalized version defined by the second author \cite{McC6, McC12}, which allow many of the results we are interested in to be stated in a slightly more streamlined fashion. Throughout, let $p$ be a prime, and let $q$ be a prime power.  Let $\mathbb{F}_q$ be the finite field with $q$ elements, and let $\widehat{\mathbb{F}_q^*}$ be the group of multiplicative characters of $\mathbb{F}_q^*$.  We extend the domain of $\chi \in \widehat{\mathbb{F}^{*}_{q}}$ to $\mathbb{F}_{q}$ by defining $\chi(0):=0$ (including for the trivial character $\varepsilon$) and denote $\bar{\chi}$ as the inverse of $\chi$. We denote by $\varphi$ the character of order two in $\widehat{\mathbb{F}_q^*}$ when $q$ is odd. More generally, for ${k>2}$ a positive integer, we let $\chi_k\in\widehat{\mathbb{F}_q^*}$ denote a character of order $k$ when $q \equiv 1 \imod{k}$. Let $\theta$ be a fixed non-trivial additive character of $\mathbb{F}_q$, and for $\chi \in \widehat{\mathbb{F}^{*}_{q}}$ define the Gauss sum $g(\chi):= \sum_{x \in \mathbb{F}_q} \chi(x) \theta(x)$.
We define the finite field hypergeometric function as follows. 
\begin{definition}[\cite{McC6}, Def. 1.4; \cite{McC12}, Def 2.4]\label{def_HypFnFF}

For $A_1, A_2, \dotsc, A_m, B_1, B_2 \dotsc, B_m \in \widehat{\mathbb{F}_q^{*}}$ and $x \in \mathbb{F}_q$,
\begin{equation*}
{_{m}F_{m}} {\biggl( \begin{array}{cccc} A_1, & A_2, & \dotsc, & A_m \\
 B_1, & B_2, & \dotsc, & B_m \end{array}
\Big| \; x \biggr)}_{q}\\
:= \frac{-1}{q-1}  \sum_{\chi \in \widehat{\mathbb{F}_q^{*}}}
\prod_{i=1}^{m} \frac{g(A_i \chi)}{g(A_i)} \frac{g(\bar{B_i \chi})}{g(\bar{B_i})}
 \chi(-1)^{m}
 \chi(x).
\end{equation*}
\end{definition}
\noindent If $B_1 = \varepsilon$, as is often the case, then it is usually omitted from the list of parameters in ${_{m}F_{m}}$ and the notation is written as ${_{m}F_{m-1}}$.
See \cite[Prop. 2.5]{McC6} and surrounding discussion for a precise description of the relationship between $_mF_m$ and Greene's function. In most of the relations connecting $_mF_m$ to Fourier coefficients of modular forms, all the $B_i$'s are trivial. In this case, $_mF_{m-1}$ equals $(-q)^{m-1}$ times Greene's function, with the same parameters. Many of the results concerning hypergeometric functions over finite fields that we quote in this paper, from other articles, were originally stated using Greene’s function. If this is the case, note then that we have reformulated them in terms of $_mF_m$, as defined above.

All the results we will see relating the $p$-th Fourier coefficients of modular forms to $_mF_m(\cdots)_p$ will require characters of certain orders. Consequently, this restricts these results to $p$ in certain congruence classes. In some cases, these results can be extended to all odd primes using a function which extends $_mF_m(\cdots)$ to the $p$-adic setting. Let $\mathbb{Z}_p$ denote the ring of $p$-adic integers, $\gfp{\cdot}$ denote Morita's $p$-adic gamma function, and $\omega$ denote the Teichm\"{u}ller character of $\mathbb{F}_p$, with $\bar{\omega}$ denoting its character inverse. 
For $x \in \mathbb{Q}$ we let  $\left\lfloor x \right\rfloor$ denote the greatest integer less than or equal to $x$ and
$\langle x \rangle$ denote the fractional part of $x$, i.e. $x- \left\lfloor x \right\rfloor$.
\begin{definition}[\cite{McC4}, Def. 2.1; \cite{McC7}, Def. 1.1]\label{def_Gp}

Let $p$ be an odd prime. For $a_1, a_2, \dotsc, a_m$, $b_1, b_2 \dotsc, b_m \in \mathbb{Q} \cap \mathbb{Z}_p$ and $x \in \mathbb{F}_p$,
\begin{multline*}
{_{m}G_{m}}
\biggl[ \begin{array}{cccc} a_1, & a_2, & \dotsc, & a_m \\
 b_1, & b_2, & \dotsc, & b_m \end{array}
\Big| \; x \; \biggr]_p
: = \frac{-1}{p-1}  \sum_{j=0}^{p-2} 
(-1)^{jm}\;
\bar{\omega}^j(x)\\
\times \prod_{i=1}^{m} 
\frac{\biggfp{\langle a_i -\frac{j}{p-1}\rangle}}{\biggfp{\langle a_i \rangle}}
\frac{\biggfp{\langle - b_i +\frac{j}{p-1}\rangle}}{\biggfp{\langle -b_i \rangle}}
(-p)^{-\lfloor{\langle a_i \rangle -\frac{j}{p-1}}\rfloor -\lfloor{\langle -b_i \rangle +\frac{j}{p-1}}\rfloor}.
\end{multline*}
\end{definition}
A ``$q$ version'' of ${_{m}G_{m}}[\cdots]$  also exists \cite[Def. 5.1]{McC7} but is not needed here.
There is a simple relationship between $_mF_m(\cdots)_p$ and ${_{m}G_{m}}[\cdots]_p$. 
\begin{lemma}[\cite{McC7}, Lemma 3.3; \cite{McC12}, Lemma 2.5]\label{lem_G_to_F}

For a fixed odd prime $p$, let $A_i, B_j \in \widehat{\mathbb{F}_p^{*}}$ be given by $\bar{\omega}^{a_i(p-1)}$ and $\bar{\omega}^{b_j(p-1)}$ respectively, where $\omega$ is the Teichm\"{u}ller character . Then
\begin{equation*}
{_{m}F_{m}} {\biggl( \begin{array}{cccc} A_1, & A_2, & \dotsc, & A_m \\
 B_1, & B_2, & \dotsc, & B_m \end{array}
\Big| \; t \biggr)}_{p}
=
{_{m}G_{m}}
\biggl[ \begin{array}{cccc} a_1, & a_2, & \dotsc, & a_m \\
 b_1, & b_2, & \dotsc, & b_m \end{array}
\Big| \; t^{-1} \; \biggr]_p.
\end{equation*}
\end{lemma}

We recall Dedekind's eta function, which will be used to describe some of the modular forms in this paper: 
$\eta(z):= q^{\frac{1}{24}} \prod_{n\geq 1} (1 -  q^n),$
where $q:= e^{2 \pi i z}.$


\section{Weight Two Newforms}\label{sec_wt2}
For each elliptic curve $E/\mathbb{Q}$, with conductor $N_E$, the modularity theorem guarantees the existence of a weight two newform of level $N_E$ whose Fourier coefficients are given by the coefficients of the Hasse--Weil $L$-function of $E$, $L(E,s)=\sum_{n\geq1}a(E, n) \, n^{-s}$.  This function is completely determined by its coefficients at the primes, $a(E,p)$, which are related to the number of rational points $N(E,p)$ on the reduction of $E$ modulo $p$ via the formula $N(E,p)=p+1-a(E,p)$. 

As mentioned in the introduction, finite field hypergeometric functions naturally lend themselves to counting points on algebraic varieties over finite fields. In particular, $N(E,p)$, and consequently $a(E,p)$, as we will see below, can be evaluated by ${}_2F_1$ finite field hypergeometric functions. Passing through the modularity theorem then results in formulas for the $p$-th Fourier coefficients of weight two newforms in terms of these ${}_2F_1$ evaluations, and we get infinitely many such connections.

The first results relating $a(E,p)$ to ${}_2F_1$ finite field hypergeometric functions were due to Koike \cite{K}. He examined various families of curves, including the Legendre family, which yields the following result.
\begin{theorem}[Koike \cite{K}]\label{thm_Koike1}

Let $\lambda\in\mathbb{Q}\setminus\{0,1\}$.  Consider the elliptic curve $E_{\lambda}:y^2=x(x-1)(x-\lambda)$ over $\mathbb{Q}$.  If $p\geq3$ is a prime with $\mathrm{ord}_p(\lambda(\lambda-1))=0$, then
\begin{equation*}
a(E_{\lambda},p) = \phi(-1) \cdot 
{_{2}F_{1}} {\biggl( \begin{array}{cc} \varphi, & \varphi \\
 & \varepsilon \end{array}
\Big| \; \lambda \biggr)}_{p}.
\end{equation*}
\end{theorem}

\begin{example}\label{eg_Legendre}:
We take $\lambda=-1$ in the above result. The curve $E_{-1}:y^2=x^3-x$ \cite[32.a3]{LMFDB} is related via the modularity theorem to the modular form
$\eta(4z)^2 \eta(8z)^2 = \sum_{n\geq1} a_1(n) \, q^n \in S_2^\text{new}\left(\Gamma_0(32)\right)$
\cite[32.2.a.a]{LMFDB}.
Combining with Theorem \ref{thm_Koike1} we get that for all odd primes
\begin{equation*}
a_1(p) = \phi(-1) \cdot 
{_{2}F_{1}} {\biggl( \begin{array}{cc} \varphi, & \varphi \\
& \varepsilon \end{array}
\Big| \; {-1} \biggr)}_{p}.
\end{equation*}
\end{example}

Fuselier \cite{F2} examined the family $E_t:y^2=4x^3-\frac{27}{1-t}x-\frac{27}{1-t}$ and expressed $a(E_t,p)$ in terms of a ${}_2F_1$, whose parameters include characters of order 12, when $p \equiv 1 \imod {12}$. Lennon \cite{L2} generalized Fuselier's result to evaluate $a(E,p)$ for any elliptic curve $E$, when the reduction of $E$ modulo $p$ is an elliptic curve over $\mathbb{F}_p$ with $j$-invariant not equal to  0 or 1728.

\begin{theorem}[Lennon \cite{L2} \S 2.2]\label{thm_Lennon1}

Let $E/\mathbb{Q}$ be an elliptic curve. Let $p \equiv 1 \imod {12}$ be a prime such that $E_p: y^2=x^3+ax+b$ is an elliptic curve over $\mathbb{F}_p$ and $j(E_p) \neq 0,1728$. Then
\begin{equation*}
a(E,p)= 
\chi_4^{\phantom{5}} \left(-\frac{a^3}{27}\right)
{_{2}F_{1}} {\biggl( \begin{array}{cc} \chi_{12}^{\phantom{5}}, & \chi_{12}^5 \\
& \varepsilon \end{array}
\Big| \; \frac{1728}{j(E_p)} \biggr)}_{p}.
\end{equation*}
\end{theorem}
Theorem \ref{thm_Lennon1} is independent of the model for $E_p$.
The results in \cite{L2} are in fact over $\mathbb{F}_q$, for $q \equiv 1 \pmod {12}$ a prime power, and hence allow calculation of $a(E,p)$ up to sign when $p \not\equiv 1 \pmod {12}$ via the relation $a(E,p)^2 = a(E,{p^2}) +2p$. 
Theorem 1.2 of \cite{McC7} extends Theorem \ref{thm_Lennon1} to the $p$-adic setting, giving a direct evaluation of $a(E,p)$ for all primes $p>3$ and resolves this sign issue.

\begin{theorem}[McCarthy \cite{McC7} Thm 1.2]\label{thm_McCarthy1}

Let $E/\mathbb{Q}$ be an elliptic curve. Let $p>3$ be a prime such that $E_p: y^2=x^3+ax+b$ is an elliptic curve over $\mathbb{F}_p$ and $j(E_p) \neq 0,1728$. Then
\begin{equation*}\label{for_main}
a(E,p)= \phi(b) \cdot p \cdot
{_{2}G_{2}}
\biggl[ \begin{array}{cc} \frac{1}{4},  & \frac{3}{4}\vspace{.05in}\\
 \frac{1}{3},  & \frac{2}{3} \end{array}
\Big| \; {1- \frac{1728}{j(E_p)}} \; \biggr]_p. 
\end{equation*}
\end{theorem}
Again, Theorem \ref{thm_McCarthy1} is independent of the model for $E_p$. 
\begin{example}:
Consider the elliptic curve $E:y^2=x^3+27x-27$ \cite[540.a2]{LMFDB}. This is related via the modularity theorem to the modular form
$ q - q^{5} - 4q^{7} + 6q^{11} - 4q^{13} - 3q^{17} - 7q^{19} - 9q^{23} + q^{25} +\cdots =\sum_{n\geq1} a_2(n) \, q^n \in S_2^\text{new}\left(\Gamma_0(540)\right)$,
\cite[540.2.a.a]{LMFDB}.
Then, by Theorem \ref{thm_McCarthy1}, we get that for all primes ${p>3}$
\begin{equation*}
a_2(p) = \left( \tfrac{-3}{p} \right) \cdot p \cdot
{_{2}G_{2}}
\biggl[ \begin{array}{cc} \frac{1}{4},  & \frac{3}{4}\vspace{.05in}\\
 \frac{1}{3},  & \frac{2}{3} \end{array}
\Big| \; {-\tfrac{1}{4}} \; \biggr]_p. 
\end{equation*}
\end{example}
The theorem covers ${p>5}$. We manually check that the relation also holds for ${p=5}$.

Significant other contributions in this area, where connections between various families of elliptic curves and finite field  hypergeometric functions are established, can be found in \cite{BK, BK2, EO2, KK, L, O}. Between all these results it should be possible to evaluate the $p$-th Fourier coefficients of all weight two newforms, with integer coefficients, using finite field hypergeometric functions (or their $p$-adic extensions). 
For newforms with non-integral coefficients, things are less straightforward and little is known. However, we have the following conjectural relations due to Evans \cite{E}.
  
\begin{conjecture}[Evans \cite{E}]\label{conj_Evans_Wt2}

Consider the newforms
$\sum_{n\geq1} a_3(n)q^n\in S_2\left(\Gamma_0(972)\right)$, \cite[972.2.a.e]{LMFDB}, and
$\sum_{n\geq1} a_4(n)q^n\in S_2\left(\Gamma_0(768)\right)$, \cite[768.2.a.j]{LMFDB},
with coefficient fields $\mathbb{Q}\left(\sqrt{2}\right)$ and $\mathbb{Q}\left(\sqrt{3}\right)$ respectively.
\begin{enumerate}
\item If $q\equiv1\pmod{6}$, then 
\begin{multline*}
\overline{\chi_6}(12)J\left(\chi_6,\chi_6\right)-\overline{\chi_6}(3) \, J\left(\chi_6^2,\chi_6^2\right)\cdot
{}_3F_2 {\biggl( \begin{array}{ccc} \overline{\chi_6}, & \varphi, & \chi_6 \\
\phantom{A_0} & \varphi \chi_6, &  \varphi \chi_6 \end{array}
\Big| \; \frac{1}{4} \biggr)}_{q}\\
=
\begin{cases}
a_3(p),&\mbox{if }q=p,\, p \equiv 1 \imod{6}\\
a_3(p)^2+2p,&\mbox{if }q=p^2,\, p \equiv 5 \imod{6}.
\end{cases}
\end{multline*}
\item If $q\equiv1\pmod{8}$, then
\begin{multline*}
\chi_8(-4)J\left(\chi_8,\chi_8\right)
-\chi_8(-4)J\left(\chi_8^2,\chi_8^3\right)\cdot{}
_3F_2 {\biggl( \begin{array}{ccc} \overline{\chi_8}, & \chi_8^3, & \chi_8 \\
\phantom{A_0} & \overline{\chi_8}^2, & \varphi \chi_8 \end{array}
\Big| \; \frac{1}{4} \biggr)}_{q}\\
=
\begin{cases}
a_4(p),&\mbox{if }q=p,\, p \equiv 1 \imod{8}\\
a_4(p)^2+2p^2,&\mbox{if }q=p^2, p \not\equiv 1 \imod{8}.
\end{cases}
\end{multline*}
\end{enumerate}
\end{conjecture}


\section{Higher Weight Newforms}\label{sec_higherwt}

The lack of universal modularity results for algebraic varieties of dimension greater than 1 means that connections to modular forms of weight greater than two are somewhat ad hoc. In this section we outline the main results linking finite field hypergeometric functions and Fourier coefficients of modular forms of weight greater than two. We start with the connections coming from Rodriguez Villegas' supercongruence conjectures.

\subsection{The Conjectures of Rodriguez Villegas}\label{subsec_RV}
In \cite{R}, Rodriguez Villegas examined the relationship between the number of points over $\mathbb{F}_p$ on certain Calabi-Yau manifolds and truncated classical hypergeometric series which correspond to a particular period of the manifold. In doing so, he identified numerically 22 possible supercongruences, 18 of which relate truncated classical hypergeometric series to Fourier coefficients of modular forms of weights three and four. 
One of these relations had previously been conjectured in \cite{SB}. The 14 cases involving weight four modular forms relate to Calabi-Yau threefolds, 13 of which were studied in \cite{BvS}. The book of Meyer \cite{Me} contains a nice description of these threefolds.

While the supercongruence relations of Rodriguez Villegas are congruences involving classical hypergeometric series, it became obvious from the work of Mortenson \cite{M} and Kilbourn \cite{Ki} on proving the first few of these conjectures, and also from known connections between finite field hypergeometric functions and some of the modular forms in question, due to Ono \cite{O} and Ahlgren and Ono \cite{AO}, that corresponding to each of Rodriguez Villegas's conjectures was a linear relation between finite field hypergeometric functions and the Fourier coefficients of these modular forms. In fact, following the work of the second author in \cite{McC4}, which precisely describes the relationship between the truncated classical hypergeometric series appearing in Rodriguez Villegas's conjectures and the $p$-adic function defined in Definition \ref{def_Gp} above, proving relationships between the relevant ${_{m}G_{m}}[\cdots]$ and the Fourier coefficients of the modular forms in question would suffice to prove the conjectures of Rodriguez Villegas. A list of all these relations is shown in Table \ref{table_RV}. The Rodriguez Villegas conjectures corresponding to cases 1-5, 7 and 15 were proved individually \cite{A, FM, I, Ki, McC5, M, VH}.  A proof of all 14 weight four conjectures (cases 5-18) is offered in \cite{LTYZ} (as yet unpublished) and, as a consequence of the results therein, all the relations in Table \ref{table_RV} should now be known. 
\newpage
\noindent
\begin{table}[htp!]
\centering
\caption{Relations arising from the conjectures of Rodriguez Villegas}
\label{table_RV}
\resizebox{\textwidth}{!}{%
\begin{tabular}{| c | c | c | c | c | c |}
\hline
&
\textbf{Hyp Series} &
\multicolumn{2}{| c }{\textbf{Newform $f(z)=\sum a(n) q^n$}} &
\multicolumn{2}{| c |}{\textbf{Connection}}\\
\hline
&
\textbf{Parameters} &
\textbf{Space} &
\textbf{LMFDB} &
\textbf{Relationship} &
\textbf{When}\\

\hline

1. &
$\left[\tfrac{1}{2}, \tfrac{1}{2}, \tfrac{1}{2} ; 1, 1, 1 | 1 \right]$ &
$S_3(\Gamma_0(16),(\tfrac{-4}{\cdot}))$ &
16.3.c.a &
$a(p) = G[\cdots]_p$ &
$p>2$
\\[6pt]

2. &
$\left[\tfrac{1}{2}, \tfrac{1}{3}, \tfrac{2}{3} ; 1, 1, 1 | 1 \right]$ &
$S_3(\Gamma_0(12),(\tfrac{-3}{\cdot}))$ &
12.3.c.a &
$a(p) = G[\cdots]_p$ &
$p>3$
\\[6pt]

3. &
$\left[\tfrac{1}{2}, \tfrac{1}{4}, \tfrac{3}{4} ; 1, 1, 1 | 1 \right]$ &
$S_3(\Gamma_0(8),(\tfrac{-2}{\cdot}))$ &
8.3.d.a &
$a(p) = G[\cdots]_p$ &
$p>2$
\\[6pt]

4. &
$\left[\tfrac{1}{2}, \tfrac{1}{6}, \tfrac{5}{6} ; 1, 1, 1 | 1 \right]$ &
$S_3(\Gamma_0(144),(\tfrac{-4}{\cdot}))$ &
144.3.g.a &
$a(p) = G[\cdots]_p$ &
$p>3$
\\[6pt]

\hline

5. &
$\left[ \tfrac{1}{2}, \tfrac{1}{2}, \tfrac{1}{2}, \tfrac{1}{2} ; 1, 1, 1, 1 | 1 \right]$ &
$S_4(\Gamma_0(8))$ &
8.4.a.a &
$a(p) = G[\cdots]_p - p$ &
$p>2$
\\[6pt]
 
6. &
$\left[ \tfrac{1}{2}, \tfrac{1}{2}, \tfrac{1}{3}, \tfrac{2}{3} ; 1, 1, 1, 1 | 1 \right]$ &
$S_4(\Gamma_0(36))$ &
36.4.a.a &
$a(p) = G[\cdots]_p - (\tfrac{12}{p}) p $ &
$p>3$
\\[6pt]
 
7. &
$\left[ \tfrac{1}{2}, \tfrac{1}{2}, \tfrac{1}{4}, \tfrac{3}{4} ; 1, 1, 1, 1 | 1 \right]$ &
$S_4(\Gamma_0(16))$ &
16.4.a.a &
$a(p) = G[\cdots]_p - (\tfrac{8}{p}) p$ &
$p>2$
\\[6pt]
 
8. &
$\left[ \tfrac{1}{2}, \tfrac{1}{2}, \tfrac{1}{6}, \tfrac{5}{6} ; 1, 1, 1, 1 | 1 \right]$ &
$S_4(\Gamma_0(72))$ &
72.4.a.b &
$a(p) = G[\cdots]_p - p$ &
$p>3$
\\[6pt]
 
 9. &
$\left[ \tfrac{1}{3}, \tfrac{2}{3}, \tfrac{1}{3}, \tfrac{2}{3} ; 1, 1, 1, 1 | 1 \right]$ &
$S_4(\Gamma_0(27))$ &
27.4.a.a&
$a(p) = G[\cdots]_p - p$ &
$p \neq 3$
\\[6pt]
 
10. &
$\left[ \tfrac{1}{3}, \tfrac{2}{3}, \tfrac{1}{4}, \tfrac{3}{4} ; 1, 1, 1, 1 | 1 \right]$ &
$S_4(\Gamma_0(9))$ &
9.4.a.a &
$a(p) = G[\cdots]_p - (\tfrac{24}{p}) p$ &
$p>3$
\\[6pt]
 
11. &
$\left[ \tfrac{1}{3}, \tfrac{2}{3}, \tfrac{1}{6}, \tfrac{5}{6} ; 1, 1, 1, 1 | 1 \right]$ &
$S_4(\Gamma_0(108))$ &
108.4.a.a &
$a(p) = G[\cdots]_p - (\tfrac{12}{p}) p$ &
$p>3$
\\[6pt]

12. &
$\left[ \tfrac{1}{4}, \tfrac{3}{4}, \tfrac{1}{4}, \tfrac{3}{4} ; 1, 1, 1, 1 | 1 \right]$ &
$S_4(\Gamma_0(32))$ &
32.4.a.a &
$a(p) = G[\cdots]_p -p$ &
$p>2$
\\[6pt]

13. &
$\left[ \tfrac{1}{4}, \tfrac{3}{4}, \tfrac{1}{6}, \tfrac{5}{6} ; 1, 1, 1, 1 | 1 \right]$ &
$S_4(\Gamma_0(144))$ &
144.4.a.f &
$a(p) = G[\cdots]_p. - (\tfrac{8}{p}) p$ &
$p>3$
\\[6pt]

14. &
$\left[ \tfrac{1}{6}, \tfrac{5}{6}, \tfrac{1}{6}, \tfrac{5}{6} ; 1, 1, 1, 1 | 1 \right]$ &
$S_4(\Gamma_0(216))$ &
216.4.a.c &
$a(p) = G[\cdots]_p - p$ &
$p>3$
\\[6pt]
 

15. & 
$\left[ \tfrac{1}{5}, \tfrac{2}{5}, \tfrac{3}{5}, \tfrac{4}{5} ; 1, 1, 1, 1 | 1 \right]$ &
$S_4(\Gamma_0(25))$ &
25.4.a.b&
$a(p) = G[\cdots]_p - (\tfrac{5}{p}) p$ &
$p \neq 5 $
\\[6pt]

16 &.
$\left[ \tfrac{1}{8}, \tfrac{3}{8}, \tfrac{5}{8}, \tfrac{7}{8} ; 1, 1, 1, 1 | 1 \right]$ &
$S_4(\Gamma_0(128))$ &
128.4.a.b &
$a(p) = G[\cdots]_p  - (\tfrac{8}{p}) p$ &
$p>2$
\\[6pt]
 
17 &.
$\left[ \tfrac{1}{10}, \tfrac{3}{10}, \tfrac{7}{10}, \tfrac{9}{10} ; 1, 1, 1, 1 | 1 \right]$ &
$S_4(\Gamma_0(200))$ &
200.4.a.f &
$a(p) = G[\cdots]_p - p$ &
$p \neq 2,5$
\\[6pt]

18. &
$\left[ \tfrac{1}{12}, \tfrac{5}{12}, \tfrac{7}{12}, \tfrac{11}{12} ; 1, 1, 1, 1 | 1 \right]$ &
$S_4(\Gamma_0(864))$ &
864.4.a.a &
$a(p) = G[\cdots]_p - p$ &
$p>3$
\\[3pt]

\hline

\end{tabular}}
\end{table}

As noted above, a couple of the relationships described in Table \ref{table_RV} were known independently of the conjectures of Rodriguez Villegas \cite{AO, O}. For example, case 5 was first proved by Ahlgren and Ono \cite{AO}.
\begin{theorem}[Ahlgren \& Ono \cite{AO}, Thm. 6]\label{thm_AO_4F3}\label{thm_AO}

Consider the weight four newform $\eta^4(2z) \, \eta^4(4z)=\sum_{n\geq1}a_5(n)q^n$ in $S_4\left(\Gamma_0(8)\right).$  If $p$ is an odd prime, then \begin{equation*}
a_5(p)=
{_{4}F_{3}} {\biggl( \begin{array}{cccc} \varphi, & \varphi, & \varphi, & \varphi \\
& \varepsilon, & \varepsilon, & \varepsilon \end{array}
\Big| \; 1 \biggr)}_{p}
-p
\end{equation*}
\end{theorem} 
Theorem \ref{thm_AO_4F3} is equivalent to case 5 in Table \ref{table_RV} via Lemma \ref{lem_G_to_F}.


\subsection{Conjectures of Evans}\label{subsec_Evans}

In addition to Conjecture \ref{conj_Evans_Wt2}, Evans also provides three conjectural relations to weight three newforms. The following is one example.

\begin{conjecture}[Evans \cite{E}]\label{conj_Evans_Wt3}

Consider the newform
$\sum_{n\geq1} a_6(n)q^n\in S_3\left(\Gamma_0(12),\left(\tfrac{-1}{\cdot}\right)\right)$, \cite[12.3.d.a]{LMFDB},
with coefficient field $\mathbb{Q}\left(\sqrt{-3}\right)$. 
If $q\equiv1\pmod{4}$, then 
$$-q-J\left(\bar{\chi_4},\bar{\chi_4}\right) \cdot
{}_3F_2 {\biggl( \begin{array}{ccc} \overline{\chi_4}, & \overline{\chi_4}, & \overline{\chi_4} \\
\phantom{\overline{\chi_4}} & \varepsilon, & \chi_4 \end{array}
\Big| \; \frac{1}{4} \biggr)}_{q}
=\begin{cases}
a_6(p),&\mbox{if }q=p,\, p\equiv1\hspace{-.3cm}\pmod{4},\\
a_6(p)^2+2p^2,&\mbox{if }q=p^2,\, p\equiv3\hspace{-.3cm}\pmod{4}.
\end{cases}$$
\end{conjecture}

The other two conjectures are similar and relate to the newforms  in $S_3\left(\Gamma_0(12),\left(\tfrac{-1}{\cdot}\right)\right)$, \cite[12.3.d.a]{LMFDB}, and
$S_3\left(\Gamma_0(972),\left(\tfrac{\cdot}{3}\right)\right)$ \cite[972.3.c.f]{LMFDB},
both with coefficient field $\mathbb{Q}\left(\sqrt{-1}\right)$.


\subsection{Relations with Ramanujan's $\tau$-function}\label{subsec_tau}

Ramanujan's $\tau$-function, $\tau(n)$, can defined as the coefficients of the unique normalized cusp form of weight 12 on the full modular group. i.e.,  $\eta(z)^{24} =: \sum_{n\geq1} \tau(n) \, q^n.$
The first result linking the $\tau$-function to finite field hypergeometric functions was given by Papanikolas.

\begin{theorem}[Papanikolas \cite{P}, Theorem 1.1]\label{thm_tau_Pap}

Let $p$ be an odd prime.  Choose $a,b\geq0$ satisfying $p=a^2+b^2$, if $p\equiv1\pmod{4}$, or $a=b=0$, if $p\equiv3\pmod{4}$.  Then 
\begin{multline*}
$$\tau(p)=-1-\left(1+\tfrac{3}{2}\phi(-1)\right)p^5+40p^3a^2b^2-128pa^4b^4\\
-\tfrac{1}{2}\sum_{\lambda=2}^{p-1}
R\left(p,\phi(1-\lambda)\cdot
{_{3}F_{2}} {\biggl( \begin{array}{ccc} \varphi, & \varphi, & \varphi \\
& \varepsilon, & \varepsilon \end{array}
\Big| \; \lambda \biggr)}_{p}
\right),
\end{multline*}
where $R(p,x)=x^5-4px^4+2p^2x^3+5p^3x^2-2p^4x-p^5$.
\end{theorem}

A similar result by Fuselier \cite{F2}, involving powers of $_{2}F_{1}(\cdots)_p$ with characters of order 12, followed. Both of these results were established using the Eichler--Selberg trace formula, which we will discuss in Section \ref{sec_Trace}.


\subsection{Other Relations}\label{subsec_Other}

In this section we mention some other noteworthy relations. Frechette, Ono and Papanikolas \cite{FOP} provide the following relation between the Fourier coefficients of a weight 6 newform and a linear combination of a ${}_6F_5$ and a ${}_4F_3$.

\begin{theorem}[Frechette, Ono \& Papanikolas \cite{FOP}, Corollary 1.2]\label{thm_FOP}

Let $\eta(z)^8\eta(4z)^4+8\eta(4z)^{12}=\sum_{n\geq1}b(n)q^n$ be the unique newform in $S_6\left(\Gamma_0(8)\right)$ \cite[6.8.a.a]{LMFDB}.  If $p$ is an odd prime, then
\begin{multline*}
b(p)={}_6F_5\biggl(\begin{array}{cccccc}\phi,&\phi,&\phi,&\phi,&\phi,&\phi\\&\varepsilon,&\varepsilon,&\varepsilon,&\varepsilon,&\varepsilon\end{array}\Big|\,1 \biggr)_p
-p\cdot{}_4F_3\biggl(\begin{array}{cccc}\phi,&\phi,&\phi,&\phi\\&\varepsilon,&\varepsilon,&\varepsilon\end{array}\Big| \,1\biggr)_p
+(1-\phi(-1)) \, p^2.
\end{multline*}
\end{theorem}
This is the only result, that we're aware of, involving a modular form of weight greater than four which can expressed via a simple linear relation of finite field hypergeometric functions.

In \cite{MP}, Papanikolas and the second author provide evidence that the eigenvalues, of index $p$, of a certain Siegel eigenform can be evaluated by the function ${}_4F_3\left(\phi, \phi, \phi, \phi ; \varepsilon, \varepsilon, \varepsilon \mid {-1}\right)_p.$
In the course of their work they prove the following.
\begin{theorem}[McCarthy \& Papanikolas \cite{MP}, Theorem 1.8]\label{thm_MP}

Consider the newform $\sum_{n\geq1} c(n) \, q^n=q+4iq^3+2q^5-8iq^7+\cdots$ in $S_3\left(\Gamma_0(32),\left(\frac{-4}{\cdot}\right)\right)$ \cite[32.3.c.a]{LMFDB}. If $p\equiv1\pmod{4}$ is prime, then 
\begin{equation*}
c(p)={}_3F_2\biggl(\begin{array}{ccc}\chi_4,&\phi,&\phi\\&\varepsilon,&\varepsilon\end{array}\Big|\,1\biggr)_p.
\end{equation*}
\end{theorem}
Both Theorems \ref{thm_FOP} and \ref{thm_MP} were proved using the Eichler--Selberg trace formula, which we will discuss in Section \ref{sec_Trace}.

In \cite{Web}, Ono provides relations for the Fourier coefficients of the only four weight three newforms which can be expressed as eta-products, all of which have complex multiplication. One such relation is as follows.
\begin{theorem}[Ono \cite{Web}, Corollary~11.20]\label{thm_Ono_CM}

Consider the newform $\eta(z)^3\eta(7z)^3 = \sum_{n\geq1} d(n) \, q^n \in S_3\left(\Gamma_0(7),\left(\frac{-7}{\cdot}\right)\right)$, \cite[7.3.b.a]{LMFDB}.
If $p\not\in\{2,3,7\}$ is prime, then
$$d(p)=\phi_p(-7) \cdot
{}_3F_2\biggl(\begin{array}{ccc}\phi,&\phi,&\phi\\&\varepsilon,&\varepsilon\end{array}\Big|\,64\biggr)_p
-\phi_p(-7) \, p.$$
\end{theorem}
The results for the other three newforms are similar. The eta-products for these newforms are $\eta(4z)^6, \eta(2z)^3\eta(6z)^3$ and $\eta(z)^2\eta(2z)\eta(4z)\eta(8z)^2$ and these are the newforms in cases 1-3 of Table \ref{table_RV} respectively.


\section{Trace Formulas for Hecke Operators}\label{sec_Trace}

There have been two main ways in which relations between finite field hypergeometric functions and Fourier coefficients of modular forms have been established. The first is via the (known or independently established) modularity of some variety, as we saw in the case of elliptic curves in Section \ref{sec_wt2}. The second is via the Eichler--Selberg trace formula for Hecke operators. These traces are connected to hypergeometric values by counting isomorphism classes of members of certain families of elliptic curves with prescribed torsion. This is a long and tedious process and works best when the dimension of the space in question is small, allowing the Fourier coefficients of specific forms to be isolated. The trace formula has also been used to establish modularity of certain varieties \cite{A2, AO2}, with the connection to hypergeometric functions following later, as was the case in Theorem \ref{thm_AO}.

For a positive integer $n$, let Tr${}_k\left(\Gamma_0(N),n\right)$ denote the trace of the $n$-th Hecke operator acting on $S_k\left(\Gamma_0(N)\right)$.
A typical result relating Tr${}_k\left(\Gamma_0(N),n\right)$ to finite field hypergeometric functions is as follows.
\begin{theorem}[Papanikolas \cite{P}, Theorem 3.2]\label{thm_Trace_Pap}

Let $p$ be an odd prime.  Choose $a,b\geq0$ satisfying $p=a^2+b^2$, if $p\equiv1\pmod{4}$, or $a=b=0$, if $p\equiv3\pmod{4}$. Define the polynomial
$$G_k(s,p)=\sum_{j=0}^{\frac{k}{2}-1}(-1)^j{{k-2-j}\choose j}p^js^{k-2j-2}.$$ 
Let 
$$\delta_k(p):=
\begin{cases}
\frac{1}{2}G_k(p,2a)+\frac{1}{2}G_k(p,2b),&\mbox{if }p\equiv1\pmod{4},\\
(-p)^{\frac{k}{2}-1},&\mbox{if }p\equiv3\pmod{4}
\end{cases}$$
and 
$$R_k(p,x):=\sum_{k=0}^{\frac{k}{2}-1}c_d\left(\tfrac{k}{2}-1\right)p^{\frac{k}{2}-1-d}x^d,$$
where $c_d(r)$ is defined by the generating function
$\frac{x+1}{\left(x^2+x+1\right)^{d+1}}=\sum_{j=-d}^\infty c_d(d+j)x^j.$
For $k\geq4$ even,
$$\mathrm{Tr}_k\left(\Gamma_0(2),p\right)=-2-\delta_k(p)-\sum_{\lambda=2}^{p-1}R_k\left(p,\phi(1-\lambda)
\cdot
{_{3}F_{2}} {\biggl( \begin{array}{ccc} \varphi, & \varphi, & \varphi \\
& \varepsilon, & \varepsilon \end{array}
\Big| \; \lambda \biggr)}_{p}
\right).$$
\end{theorem}
Taking $k=12$ in Theorem \ref{thm_Trace_Pap}, and using the fact that $\eta(z)^{24}$ and $\eta(2z)^{24}$ form a basis for $S_{12}\left(\Gamma_0(2)\right)$, yields Theorem \ref{thm_tau_Pap}. By taking $k=8$ and $k=10$ in Theorem \ref{thm_Trace_Pap}, Papanikolas also provides formulas, similar to that in Theorem \ref{thm_tau_Pap}, for the coefficients of the unique newforms in $S_{8}\left(\Gamma_0(2)\right)$ and $S_{10}\left(\Gamma_0(2)\right)$ respectively.

Similar evaluations of the traces of the $p$-th Hecke operators acting on the following spaces have also been produced.
\begin{itemize}

\item $S_{k}\left(\Gamma_0(4)\right)$, for $k\geq 4$ even \cite{A2, FOP};\\[-4pt]

\item $S_{k}\left(\Gamma_0(8)\right)$, for $k\geq 4$ even \cite{AO2, FOP};\\[-4pt]

\item $S_{k}\left(\Gamma \right)$, for $k\geq 4$ even \cite{F2, F3};\\[-4pt]

\item $S_{k}\left(\Gamma_0(3)\right)$, $S_{k}\left(\Gamma_0(9)\right)$, for $k\geq 4$ even \cite{L}; and\\[-4pt]

\item $S_3\left(\Gamma_0(16),\left(\tfrac{-4}{\cdot}\right)\right)$, $S_3\left(\Gamma_0(32),\left(\tfrac{-4}{\cdot}\right)\right)$ \cite{MP}.

\end{itemize}
Lennon \cite{L} uses the evaluation for Tr${}_k\left(\Gamma_0(9),p\right)$, when ${k=4}$, to give another formula for the $p$-th Fourier coefficients of the newform in case 10 of Table \ref{table_RV}.

\begin{theorem}[Lennon \cite{L}, Corollary 1.8]\label{thm_Lennon2}

Let $\eta(3z)^8=\sum_{n\geq1} h(n) \, q^n \in S_4^\text{new}\left(\Gamma_0(9)\right)$, \cite[9.4.a.a]{LMFDB}. For $p\equiv1\pmod{3}$,
$$h(p)={}_2F_1\biggl(\begin{array}{cc}\chi_3,&\bar{\chi_3}\\&\varepsilon\end{array}\Big|\,9\cdot8^{-1}\biggr)_{p^3}.$$
\end{theorem}


\section{New Relations}\label{sec_NewConjs}
In recent work \cite{DM}, we examined the number of complete subgraphs of order four contained in generalized Paley graphs. Let $k \geq 2$ be an integer. Let $q$ be a prime power such that $q \equiv 1 \imod {k}$ if $q$ is even, or, $q \equiv 1 \imod {2k}$ if $q$ is odd. The generalized Paley graph of order $q$, $G_k(q)$, is the graph with vertex set $\mathbb{F}_q$ where $ab$ is an edge if and only if ${a-b}$ is a $k$-th power residue. We provided a formula, in terms of ${}_3F_2$ finite field hypergeometric functions, for the number of complete subgraphs of order four contained in $G_k(q)$, which holds for all $k$. This formula includes all
$${_{3}F_2}\biggl( \begin{array}{ccc} \chi_k^{t_1}, & \chi_k^{t_2}, & \chi_k^{t_3} \vspace{.05in}\\
\phantom{\chi_k^{t_1}} & \chi_k^{t_4}, & \chi_k^{t_5} \end{array}
\Big| \; 1 \biggr)_{q},$$
as $(t_1,t_2,t_3,t_4,t_5)$ ranges over all tuples in $\left( \mathbb{Z} / k\mathbb{Z} \right)^{5}$.
We also showed that many of these terms can be simplified and many are equal to each other. We gave explicit determinations for $k \leq 4$ and noticed that many of the ${}_3F_2$'s that remained were known to be related to Fourier coefficients of weight three modular forms. We also found numerically two new possible relations. Specifically, consider the newform $g_1(z)  = q + 3 i q^{2} -5 q^{4} -3 i q^{5} + 5 q^{7} -3 i q^{8} + \cdots = \sum_{n=1}^{\infty} \beta_1(n) q^n \in S_3(\Gamma_0(27), (\tfrac{-3}{\cdot}))$ \cite[27.3.b.b]{LMFDB}.
Then numerical evidence suggests that, for $p \equiv 1 \imod6$,
\begin{equation}\label{eqn_Paley_1}
{_{3}F_2} {\biggl( \begin{array}{cccc} \chi_3, & \chi_3, &  \bar{\chi_3} \\
\phantom{\chi_3} & \varepsilon, &  \varepsilon \end{array}
\Big| \; 1 \biggr)}_{p} 
= \beta_1(p).
\end{equation}
Also, consider the newform $g_2(z)  = q + ( 2 \zeta_{8} - 2 \zeta_{8}^{3} ) q^{3} + 4 \zeta_{8}^{2} q^{5} + ( 8 \zeta_{8} + 8 \zeta_{8}^{3} ) q^{7} - q^{9} +  \cdots = \sum_{n=1}^{\infty} \beta_2(n) q^n \in S_3(\Gamma_0(128), (\tfrac{-8}{\cdot}))$ \cite[128.3.d.c]{LMFDB} ,
for a primitive eighth root of unity $\zeta_{8}$. Then, for $p \equiv 1 \imod 4$, we observed
\begin{equation}\label{eqn_Paley_2}
{_{3}F_2} \biggl( \begin{array}{ccc} \chi_4, & \chi_4 & \bar{\chi_4} \vspace{.02in}\\
\phantom{\chi_4} & \varepsilon, & \varepsilon \end{array}
\Big| \; 1 \biggr)_{p}
=\pm  \beta_2(p).
\end{equation}
Since then we have carried out a more extensive search. We examined the ${}_3F_2$'s coming from the results in \cite{DM} for all $k\leq12$. We focused our search based on what appeared to be desirable characteristics that we observed in the small $k$ cases. The ${}_3F_2$'s can be sorted into orbits (see \cite{DM} for precise details) and all the new relations we found were where the ${}_3F_2$ and its conjugate were in the same orbit. 
All the new conjectural relations we have found are summarized in Table \ref{table_NewConjs}. To simplify the table we have listed the parameters of the ${}_3F_2$'s using rational numbers according to the convention that the fraction $\frac{t}{k}$ represents the character $\chi_k^{t}$. Interestingly, our search yielded the relations in cases 1-4 of Table \ref{table_RV} and the relation in Theorem \ref{thm_MP}, but as they are already known, we have not included them in Table \ref{table_NewConjs}. However, for completeness, we have included the relations from (\ref{eqn_Paley_1}) and  (\ref{eqn_Paley_2}). They appear as cases 4 and 8 respectively.
Similar to the relation in (\ref{eqn_Paley_2}), many of the new relations involve a sign which doesn't seem to be resolvable by a simple twist by a Dirichlet character. So we first define some functions to explain these signs.

For $p \equiv 1 \pmod{4}$, write $p=x^2+y^2$ for integers $x$ and $y$, such that $x$ is odd and $y$ is even.
For $p\equiv 1 \pmod{12}$, note that either $3 \mid x$ or $3 \mid y$ and define
\begin{equation*}
S_x(p)
=
\begin{cases}
+1 & \textup{if } 3 \mid y;\\
-1 & \textup{if } 3\mid x.
\end{cases}
\end{equation*}
Note that $S_x(p)$ equals $c_{12}^2$, where $c_{12}$ is the quantity described in \cite[Ch. 3.5]{BEW}.
If $p\equiv 1 \pmod{20}$, then  either $5 \mid x$ or $5 \mid y$ and so we define
\begin{equation*}
S_{20}(p)
=
\begin{cases}
+1 & \textup{if } 5 \mid y,\\
-1 & \textup{if } 5 \mid x.
\end{cases}
\end{equation*}
Now, for $p\equiv 1 \pmod{6}$, define
\begin{equation*}
S_6(p)
=
\begin{cases}
S_x(p) & \textup{if } p\equiv 1 \pmod{12},\\
\pm1 & \textup{if } p\equiv 7 \pmod{12}.
\end{cases}
\end{equation*}

For $p \equiv 1 \pmod{8}$, write $p=u^2+2v^2$ for integers $u$ and $v$, such that $u\equiv 3 \imod 4$ and $v$ is even. For $p\equiv 1 \pmod{4}$, define
\begin{equation*}
S_4(p)
=
\begin{cases}
+1 & \textup{if } p \equiv 1 \imod{8} \text{ and }  v \equiv 0 \imod{4}, \text{ or, } p \equiv 13 \imod{16},\\
-1 & \textup{if } p \equiv 1 \imod{8} \text{ and }  v \equiv 2 \imod{4}, \text{ or, } p \equiv 5 \imod{16}.
\end{cases}
\end{equation*}
When $p\equiv 1 \pmod{12}$, define
\begin{equation*}
S_u(p)
=
\begin{cases}
+1 & \textup{if } u \equiv 2\imod{3},\\
-1 & \textup{if } u \equiv 1\imod{3},
\end{cases}
\quad \text{and,} \quad
S_{12}(p)
=
\begin{cases}
S_u(p) & \textup{if } p\equiv 1 \pmod{24},\\
\pm1 & \textup{if } p\equiv 13 \pmod{24}.
\end{cases}
\end{equation*}

For $p \equiv 1 \pmod{10}$, write $p=a^2 + 5 b^2 + 5 c^2 + 5 d^2$ for integers $a,b,c,d$ such that $a \equiv 4 \imod 5 $ and $ab=d^2-c^2-cd$. We note that $a$ is unique up to sign \cite[Thm. 3.7.2]{BEW}. Define
\begin{equation*}
S_{10}(p)
=
\begin{cases}
+1 & \textup{if } 4 \nmid a;\\
-1 & \textup{if }  4 \mid a.
\end{cases}
\end{equation*}
$S_{10}(p)$ relates to case 14 in Table \ref{table_NewConjs}. The modular form in that case has CM by $\mathbb{Q}\left(\sqrt{-5}\right)$ and so its Fourier coefficients $a(p)$ vanish when {$p\equiv 11 \imod{20}$}. Thus we only need the sign at $p\equiv 1 \imod{20}$. When $p\equiv 1 \pmod{20}$ it appears $S_{10}(p)=S_{20}(p)$.
\begin{table}[htp!]
\caption{New Conjectural Relations}
\label{table_NewConjs}
\resizebox{\textwidth}{!}{%
\begin{tabular}{| c | c | c | c | c | c |}
\hline
&
\textbf{Hyp Series} &
\multicolumn{2}{| c }{\textbf{Newform $f(z)=\sum a(n) q^n$}} &
\multicolumn{2}{| c |}{\textbf{Connection}}\\
\hline
&
\textbf{Parameters} &
\textbf{Space} &
\textbf{LMFDB} &
\textbf{Relationship} &
\textbf{Conditions}\\

\hline

1. &
$\left[\tfrac{1}{3}, \tfrac{1}{2}, \tfrac{1}{2} ; 1, 1, 1 | 1 \right]$ &
$S_3(\Gamma_0(48), (\tfrac{-4}{\cdot}))$ &
48.3.g.a &
$a(p) = S_{6}(p) \cdot F(\cdots)_p$ &
$p \equiv 1 \imod 6$
\\[6pt]

2. &
$\left[\tfrac{1}{6}, \tfrac{1}{2}, \tfrac{1}{2} ; 1, 1, 1 | 1 \right]$ &
$S_3(\Gamma_0(12), (\tfrac{-4}{\cdot}))$ &
12.3.d.a &
$a(p) = S_{6}(p) \cdot F(\cdots)_p$ &
$p \equiv 1 \imod 6$
\\[6pt]

3. &
$\left[\tfrac{1}{8}, \tfrac{1}{2}, \tfrac{1}{2} ; 1, 1, 1 | 1 \right]$ &
$S_3(\Gamma_0(64), (\tfrac{-8}{\cdot}))$ &
64.3.d.a &
$a(p) = F(\cdots)_p$ &
$p \equiv 1 \imod 8$
\\[6pt]


4. &
$\left[\tfrac{1}{3}, \tfrac{1}{3}, \tfrac{2}{3} ; 1, 1, 1 | 1 \right]$ &
$S_3(\Gamma_0(27), (\tfrac{-3}{\cdot}))$ &
27.3.b.b &
$a(p) = F(\cdots)_p$ &
$p\equiv 1 \imod 6$
\\[6pt]

5. &
$\left[\tfrac{1}{4}, \tfrac{1}{3}, \tfrac{2}{3} ; 1, 1, 1 | 1 \right]$ &
$S_3(\Gamma_0(36), (\tfrac{-4}{\cdot}))$ &
36.3.d.a &
$a(p) = F(\cdots)_p$ &
$p\equiv 1 \imod {12}$
\\[6pt]

6. &
$\left[\tfrac{1}{6}, \tfrac{1}{3}, \tfrac{2}{3} ; 1, 1, 1 | 1 \right]$ &
$S_3(\Gamma_0(108),(\tfrac{-3}{\cdot}))$ &
108.3.c.b &
$a(p) = F(\cdots)_p$ &
$p\equiv 1 \pmod 6$
\\[6pt]


7. &
$\left[\tfrac{1}{3}, \tfrac{1}{4}, \tfrac{3}{4} ; 1, 1, 1 | 1 \right]$ &
$S_3(\Gamma_0(576),(\tfrac{-24}{\cdot}))$ &
576.3.h.b &
$a(p) = S_{12}(p) \cdot F(\cdots)_p$ &
$p \equiv 1 \imod {12}$
\\[6pt]

8. &
$\left[\tfrac{1}{4}, \tfrac{1}{4}, \tfrac{3}{4} ; 1, 1, 1 | 1 \right]$ &
$S_3(\Gamma_0(128),(\tfrac{-8}{\cdot}))$ &
128.3.d.c &
$a(p) = S_4(p) \cdot F(\cdots)_p$ &
$p \equiv 1 \imod 4$
\\[6pt]

9. &
$\left[\tfrac{1}{6}, \tfrac{1}{4}, \tfrac{3}{4} ; 1, 1, 1 | 1 \right]$ &
$S_3(\Gamma_0(576),(\tfrac{-24}{\cdot}))$ &
576.3.h.a &
$a(p) = S_{12}(p) \cdot F(\cdots)_p$ &
$p \equiv 1 \imod {12}$
\\[6pt]


10. &
$\left[\tfrac{1}{3}, \tfrac{1}{6}, \tfrac{5}{6} ; 1, 1, 1 | 1 \right]$ &
$S_3(\Gamma_0(432),(\tfrac{-4}{\cdot}))$ &
432.3.g.a &
$a(p) = S_{6}(p) \cdot F(\cdots)_p$ &
$p \equiv 1 \imod 6$
\\[6pt]

11. &
$\left[\tfrac{1}{4}, \tfrac{1}{6}, \tfrac{5}{6} ; 1, 1, 1 | 1 \right]$ &
$S_3(\Gamma_0(288),(\tfrac{-4}{\cdot}))$ &
288.3.g.a &
$a(p) = F(\cdots)_p$ &
$p \equiv 1 \imod {12}$
\\[6pt]

12. &
$\left[\tfrac{1}{6}, \tfrac{1}{6}, \tfrac{5}{6} ; 1, 1, 1 | 1 \right]$ &
$S_3(\Gamma_0(108),(\tfrac{-4}{\cdot}))$ &
108.3.d.a &
$a(p) = S_{6}(p) \cdot F(\cdots)_p$ &
$p \equiv 1 \imod 6$
\\[6pt]


13. &
$\left[\tfrac{1}{5}, \tfrac{1}{5}, \tfrac{4}{5} ; 1, 1, 1 | 1 \right]$ &
$S_3(\Gamma_0(25),\chi)$ &
25.3.c.a &
$a(p) = F(\cdots)_p$ &
$p \equiv 1 \imod 5$
\\[6pt]

14. &
$\left[\tfrac{1}{2}, \tfrac{1}{10}, \tfrac{9}{10} ; 1, 1, 1 | 1 \right]$ &
$S_3(\Gamma_0(20), (\tfrac{-20}{\cdot}))$ &
20.3.d.a &
$a(p) = S_{10}(p) \cdot F(\cdots)_p$ &
$p \equiv 1 \imod{10}$
\\[6pt]

15. &
$\left[\tfrac{1}{2}, \tfrac{1}{12}, \tfrac{11}{12} ; 1, 1, 1 | 1 \right]$ &
$S_3(\Gamma_0(24), (\tfrac{-24}{\cdot}))$ &
24.3.h.a &
$a(p) = F(\cdots)_p$ &
$p \equiv 1 \imod{12}$
\\[6pt]

\hline

\end{tabular}}
Note: $\chi$ on row 13 is the Dirichlet character of conductor 5, with $2 \mapsto i$.
\end{table}

As we have seen, the functions $S_6$ and $S_{12}$, which affect cases 1, 2, 7, 9, 10 and 12 in Table \ref{table_NewConjs}, are not fully described when {$p \equiv 7 \imod{12}$} and {$p \equiv 13 \imod{24}$} respectively. Unfortunately, we were unable to ascribe a simple formula to the sign in those classes for those cases. Also, the choice of character is important in those cases. This is best explained using case 1 as an example. Combining \cite[(4.25)]{G} with \cite[Prop. 2.5]{McC6} we see that, for $p \equiv 1 \imod{6}$,
\begin{equation*}
{_{3}F_2}\biggl( \begin{array}{ccc} \chi_3, & \varphi, & \varphi \vspace{.05in}\\
\phantom{\varepsilon} & \varepsilon, & \varepsilon \end{array}
\Big| \; 1 \biggr)_{p}
=
{_{3}F_2}\biggl( \begin{array}{ccc} \bar{\chi_3}, & \varphi, & \varphi \vspace{.05in}\\
\phantom{\varepsilon} & \varepsilon, & \varepsilon \end{array}
\Big| \; 1 \biggr)_{p}
\times
\begin{cases}
+1 & \text{if } p \equiv 1 \pmod 4,\\
-1 & \text{if } p \equiv 3 \pmod 4.
\end{cases}
\end{equation*}
There are two characters of order three when $p \equiv 1 \imod{6}$ and they are conjugates of each other. So, when $p \equiv 1 \imod{12}$, the ${_{3}F_2}$ in case 1 is independent of the choice of $\chi_3$. However, when $p \equiv 7 \imod{12}$, the choice of character will determine the sign. Similar behavior is observed in cases 2, 7, 9, 10 and 12.

It doesn't appear that the relations in Table \ref{table_NewConjs} can be extended to all primes in a simple way using ${_{m}G_{m}}[\cdots]_p$. It may be possible, however, with the introduction of extra factors which equal $\pm1$ when $p$ is in the equivalence class outlined in the table.


\section*{Acknowledgements}
The first author is supported by an AMS-Simons travel grant from the American Mathematical Society and the Simons Foundation. The second author is supported by a grant from the Simons Foundation (\#353329, Dermot McCarthy).



\end{document}